\RequirePackage{ifpdf}
\ifpdf % We are running pdfTeX in pdf mode
\documentclass[pdftex]{sigma}
\else
\documentclass{sigma}
\fi

\begin{document}
\allowdisplaybreaks

\renewcommand{\PaperNumber}{048}

\FirstPageHeading

\ShortArticleName{On Deformations and Contractions of Lie Algebras}

\ArticleName{On Deformations and Contractions of Lie Algebras}

\Author{Alice FIALOWSKI~$^\dag$ and Marc DE MONTIGNY~$^\ddag$}

\AuthorNameForHeading{A. Fialowski and M. de Montigny}

\Address{$^\dag$~Institute of Mathematics, E\"otv\"os Lor\'and University,
 P\'azm\'any P\'eter s\'et\'any 1/C, \\
 $\phantom{^\dag}$~H-1117, Budapest, Hungary} 
\EmailD{\href{mailto:fialowsk@cs.elte.hu}{fialowsk@cs.elte.hu}} 
\URLaddressD{\url{http://www.cs.elte.hu/~fialowsk/}} 

\Address{$^\ddag$~Campus Saint-Jean and Theoretical Physics Institute,
University of Alberta, 8406 - 91 Street,\\
$\phantom{^\ddag}$~Edmonton, Alberta, T6C 4G9, Canada}
\EmailD{\href{mailto:montigny@phys.ualberta.ca}{montigny@phys.ualberta.ca}} 
\URLaddressD{\url{http://www.phys.ualberta.ca/~montigny/}}

\ArticleDates{Received February 24, 2006, in f\/inal form April
25, 2006; Published online May 03, 2006}

\Abstract{In this contributed presentation, we 
 discuss and compare the mutually 
 opposite procedures of deformations and contractions
 of Lie algebras. We suggest  
 that with approp\-ria\-te combinations of both procedures one may
 construct new Lie algebras. We f\/irst discuss low-dimensional Lie
 algebras and illustrate thereby that whereas for every contraction
 there exists a reverse deformation, the converse is not true
 in general. Also we note that some Lie algebras belonging to 
 parameterized families are singled out by the irreversibility of 
 deformations and contractions. After reminding that global
 deformations of the Witt, Virasoro, and af\/f\/ine Kac--Moody algebras
 allow one to retrieve Lie algebras of Krichever--Novikov type, we
 contract the latter to f\/ind new inf\/inite dimensional Lie algebras.}


\Classification{17B66; 17B67; 17B65; 17B56; 17B68; 14D15; 81T40}

\section{Introduction}

The purpose of this presentation is to report on a recent 
 analysis of two dif\/ferent approaches to `deformations' 
 of Lie algebras, by which loose term 
 we mean continuous modif\/ications of their structure
 constants. These deformations appear in mathematics and in 
 physics under various guises. We analyze hereafter
 two main categories of such modif\/ications: {\it contractions}, which
 typically transform a Lie algebra into a `more Abelian' Lie algebra, and
 {\it deformations}, which, understood in a strict sense, lead to Lie
 algebras with more intricate Lie brackets. 
 There exists a plethora of def\/initions for both contractions
 and deformations and below, we recall general def\/initions
 within the framework of Lie algebras. In order to clarify
 the two concepts and provide explicit constructions, we
 consider concrete examples in low dimensions and inf\/inite
 dimensions. For further details
 and references, see~\cite{fialowskimontigny}.

References about early works on deformations of mathematical structures
 in general, and Lie algebras in particular, are given in~\cite{fialowskimontigny}. 
 We have listed some of them in~\cite{gerstenhaber}. Of particular interest hereafter is 
 the more general deformation theory obtained
 by considering an arbitrary commutative algebra with unity as the base
 of deformation. Such deformations are called `global' and appear in the 
 work of Fialowski and Schlichenmaier~\cite{fialowskischli}. 

A contraction is a procedure somewhat opposite to deformation. Contractions are
 important in physics because they explain in terms of Lie
 algebras why some theories arise as a limit regime of more `exact'
 theories. They consist in multiplying the generators
 of the symmetry by `contraction parameters' such that
 when these parameters reach some singularity point one obtains a 
 non-isomorphic Lie algebra with the same dimension~\cite{iw}.
 The mathematics literature contains various concepts similar to
 contractions: `degeneration', `orbit closure', etc. 
 
Now let us def\/ine both concepts.
 A review of the concepts of deformations and contractions is
 given in~\cite{vinberg}. Some articles which address various
 aspects of both deformations and contractions are in~\cite{weimar2000}.
 Consider a Lie algebra ${\mathfrak{g}}$ of dimension $N$ over the 
 f\/ield $k$ which we take hereafter as being $k={\mathbb R}$ and ${\mathbb C}$.
 We write the basis elements of ${\mathfrak{g}}$ as 
 $\{x_1,\ldots,x_N\}$ with Lie brackets
\begin{gather}
[x_i,x_j]=C_{ij}^kx_k,
\label{liebrackets}
\end{gather}
where the coef\/f\/icients $C_{ij}^k$ are the structure constants.
 We denote by ${\cal L}_N(k)$ the space of structural tensors
 of $N$-dimensional Lie algebras. A {\it one-parameter deformation} of
 a Lie algebra ${\mathfrak{g}}$, with structure constants belonging
 to ${\cal L}_N(k)$, can be seen as
 a continuous curve over ${\cal L}_N(k)$. One refers to the deformation
 as being (piecewise) smooth, analytic, etc.\ if the associated
 def\/ining curve itself is (piecewise) smooth, analytic, respectively. 

A {\it formal} one-parameter deformation is def\/ined by the Lie brackets:
\begin{gather}
[a,b]_t=F_0(a,b)+tF_1(a,b)+\cdots+t^mF_m(a,b)+\cdots,
\label{formaldeformation}
\end{gather}
 where $F_0$ denotes the original Lie bracket $[\cdot,\cdot]$
 and $F_m$ are two-cochains. The Jacobi identity implies, among others,
 that $F_1$ must be a two-cocycle of ${\mathfrak{g}}$. 
 We call $[\cdot,\cdot]_t$ a f\/irst-order, or {\it infinitesimal},
 deformation  if it satisf\/ies the Jacobi identity up to $t^2$. 
 It follows that f\/irst-order deformations correspond to elements of
 the space of two-cocycles  $Z^2({\mathfrak{g}},{\mathfrak{g}})$. 
 When a formal deformation is such that $[\;,\;]_t\simeq [\;,\;]_s$ 
 for every $t$ and $s$ except $0$, then we call it a {\it jump
 deformation}. It is clear from the def\/inition that
 contractions (def\/ined in equation~(\ref{contraction})) are 
 exactly related to jump deformations. 

Now let us look at a deformation ${\mathfrak g}_t = [\cdot,\cdot]_t$ not as a
 one-parameter 
 family of Lie algebras, but as a Lie algebra over the ring $k[[t]]$
 of formal power series over $k$.
 A natural generalization is to allow more parameters, which amounts
 to consider $k[[t_1,\dots,t_k]]$ as the base, or, even more generally,
 to take an arbitrary commutative algebra $A$ over $k$, with unit as the base.
 Assume that $A$ admits an augmentation $\epsilon : A\rightarrow k$,
 such that $\epsilon$ is a $k$-algebra homomorphism and $\epsilon(1_A) = 1$. The
 ideal $m_\epsilon :=\rm{ker}(\epsilon)$ is a maximal ideal of $A$, and, given
 a maximal ideal $m$ of $A$ with $A/m\cong k$, the natural quotient map
 def\/ines an augmentation. If $A$ has a unique maximal ideal, the deformation
 with base $A$ is called {\it local}. If $A$ is the projective limit
 of local algebras, the deformation is called {\it formal}.
 In Section~3, we will consider inf\/inite dimensional Lie algebras
 obtained by using the concept of {\it global} deformation. For more details, 
 see~\cite{fialowskimontigny}.

Intuitively, {\it rigidity} of a Lie algebra $\mathfrak g$ means that we
 cannot deform it.  We call a Lie algeb\-ra
 {\it infinitesimally rigid} if every inf\/initesimal deformation is
 equivalent to the trivial one, and {\it formally rigid} if every
 formal deformation is trivial. The examples discussed below
 are formally rigid (see~\cite{fialowski1990,formalrigid}). As mentioned in Section 3,
 the interesting feature of inf\/inite dimensional Lie algebras is that formal
 deformations are no longer suf\/f\/icient to describe general deformations. The
 examples discussed below are formally rigid
 (see~\cite{fialowski1990}), so that they admit no
 non-trivial formal deformations. Nevertheless, there exist very interesting
 non-trivial global deformations. In global deformation theory, we no longer have the tool of computing cohomology in order to get
 deformations so that the picture is much more dif\/f\/icult and there are very few
 results so far~\cite{fialowskischli}.

 The commutation relations of a {\it contracted Lie algebra},
 or {\it contraction}, ${\mathfrak{g}}'$ of a Lie algebra ${\mathfrak{g}}$, are given by 
 the limit \cite{iw,conatser,burdesteinhoff}:
\begin{gather}
[x,y]'\equiv\lim_{\varepsilon\rightarrow\varepsilon_0}
 {\cal U}_\varepsilon^{-1}([{\cal U}_\varepsilon(x),
 {\cal U}_\varepsilon(y)],
\label{contraction}
\end{gather}
where ${\cal U}_\varepsilon\in\ $GL($N,k$) is a non-singular linear transformation
 of ${\mathfrak{g}}$, with $\varepsilon_0$ being a
 singularity point of its inverse ${\cal U}_\varepsilon^{-1}$. 

Throughout the paper, however, we shall utilize 
 contractions def\/ined with {\it diagonal} ${\cal U}_\varepsilon$
 in equation~(\ref{contraction}),
 def\/ined by splitting the Lie
 algebra ${\mathfrak{g}}$ into an arbitrary number of subspaces:
\begin{gather}
{\mathfrak{g}}= {\mathfrak{g}}_0 + {\mathfrak{g}}_1 +\cdots +{\mathfrak{g}}_p,\label{wwsplit}
\end{gather}
 and by taking the matrix ${\cal U}_\varepsilon$ as follows:
 \begin{gather}
{\cal U}_\varepsilon^{\rm diag}=
\oplus_j\ \varepsilon^{n_j}\ {\mathrm {id}}_{\mathfrak{g}_j},\qquad
  \varepsilon>0, \quad n_j\in {\mathbb R} ,\quad j=1,2,\dots, p,\label{weimarwoods}
 \end{gather}
 where $p\leq$dim$\; {\mathfrak{g}}$. 
 From equations (\ref{liebrackets}) and (\ref{weimarwoods}), and if we
 denote by $\mathfrak g_i$ the subspace in (\ref{wwsplit}) to which the element
 $x_i$ belongs, then equation~(\ref{contraction}) becomes
\begin{gather*}
[x_i,x_j]'=\lim_{\varepsilon\rightarrow 0} \varepsilon^{n_i+n_j-n_k} C^k_{ij} x_k.
%\label{wwpowers}
\end{gather*}
Thus the exponents in equation~(\ref{weimarwoods}) must satisfy 
\begin{gather}
n_i+n_j-n_k\geq 0,
\label{powerrelation}
\end{gather}
 unless $C^k_{ij}=0$.
 Then the structure constants of the contracted algebra
 $\mathfrak g'$ are given by
 $(C')^k_{ij}= C^k_{ij}$ if $n_i+n_j=n_k$, and
 $(C')^k_{ij}=0$ if $n_i+n_j>n_k$.
Two trivial contractions are always present: the Abelian Lie algebra
 and the original Lie algebra itself, for which the commutation relations are
 unchanged. Likewise, an Abelian Lie algebra can be deformed to every
 Lie algebra of the same dimension.

\section{Three-dimensional complex Lie algebras}

 In this section, we enumerate the deformations and the contractions
 of complex three-dimensional Lie algebras, in order to demonstrate
 the dif\/ferences between the two concepts. 
 Deformation of complex three-dimensional
 Lie algebras were recently classif\/ied in~\cite{fialowskipenkava}.  
 The real algebras are discussed in~\cite{fialowskimontigny}.  
 
 \begin{table}
\caption{Three-dimensional complex Lie algebras.}\label{3dcomplex}
\vspace{1mm}

\centering\begin{tabular}{|ll|}\hline
 & \\
${\mathbb C}^{3}:$ & $[x_{i},x_{j}]=0,\qquad i,j=1,2,3$\\ 
$\mathfrak{n}_{3}\left( {\mathbb C}\right): $ & $[x_{1},x_{2}]=x_{3}$ \\ 
$\mathfrak{r}_{2}\left( {\mathbb C}\right) \oplus {\mathbb C}:$ & $[x_{1},x_{2}]=x_{2}$ \\ 
$\mathfrak{r}_{3}\left( {\mathbb C}\right): $ & $%
[x_{1},x_{2}]=x_{2},[x_{1},x_{3}]=x_{2}+x_{3}$ \\ 
$\mathfrak{r}_{3,\lambda }\left( {\mathbb C}\right), (\lambda\in {\mathbb C}^*,
 |\lambda|\leq 1): $ & $%
[x_{1},x_{2}]=x_{2},[x_{1},x_{3}]=\lambda x_{3}$ \\
$\mathfrak{sl}_{2}\left( {\mathbb C}\right): $ & $%
[x_{1},x_{2}]=x_{3},[x_{2},x_{3}]=x_{1},[x_{3},x_{1}]=x_{2}$\\
 &  \\
 \hline
\end{tabular}

\end{table}

 The Lie brackets of the three-dimensional complex Lie algebras are given in 
 Table~\ref{3dcomplex}. Note that ${\mathfrak {r}}_{3,{\bar{\lambda}}}({\mathbb C})$ is
 isomorphic to ${\mathfrak {r}}_{3,{{\lambda}}}({\mathbb C})$ when $|\lambda|=1$
 because $\lambda{\bar{\lambda}}=1$. 
 As a simple illustration of the methods, consider the
 contraction from ${\mathfrak{sl}_2}({\mathbb C})$  to ${\mathfrak {r}}_{3,-1}({\mathbb C})$. We express 
 the Lie brackets of ${\mathfrak{sl}_2}({\mathbb C})$ in the Cartan basis:
\begin{gather}
[h,e]=e,\qquad [h,f]=-f,\qquad [e,f]=2h.
\label{sl2efh}
\end{gather} 
Then, we may introduce the contraction parameters as follows: 
\[
e\rightarrow\varepsilon e,\qquad
f\rightarrow\varepsilon f,\qquad
h\rightarrow h,
\]
before taking the limit $\varepsilon\rightarrow 0$. This results in
 $[e,f]\rightarrow 0$, with $[h,e]$ and $[h,f]$ unchanged,
 i.e.\ the Lie brackets for ${\mathfrak {r}}_{3,-1}({\mathbb C})$. 
 Now, let us illustrate the reverse deformation with this simple example. 
 The original Lie brackets of ${\mathfrak {r}}_{3,-1}({\mathbb C})$ are such that, in
  equation (\ref{formaldeformation}), the non-zero $F_0$'s are 
 $F_0(h,e)=e$  and $F_0(h,f)=-f$. Then, in order to deform it
 to ${\mathfrak{sl}_2}({\mathbb C})$, we may write equation~(\ref{formaldeformation}) as
\[
[h,e]_t=e+tF_1(h,e),\qquad
[h,f]_t=-f+tF_1(h,f),\qquad
[e,f]_t=tF_1(e,f),
\]
where
\[
F_1(h,e)=0,\qquad F_1(h,f)=0,\qquad F_1(e,f)=2h,
\]
as suggested clearly by the contraction. The resulting Lie algebra is
 isomorphic to ${\mathfrak{sl}_2}({\mathbb C})$, for any non-zero $t$.

\begin{figure}[t]
\centerline{\includegraphics{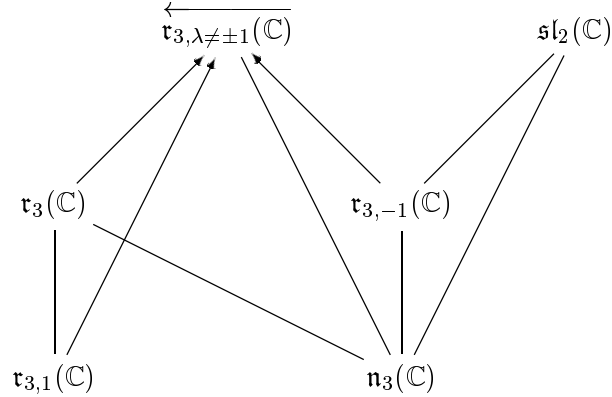}}
\caption{Contractions and deformations of
 the three-dimensional complex Lie algebras.}\label{fig1}
 \end{figure}

The results of contractions and deformations of three-dimensional complex Lie
 algebras are displayed on Fig.~\ref{fig1}.
 The lines and arrows should be interpreted as follows:
 an arrow points toward the deformation, whereas a simple line connects Lie algebras
 related by both deformation and contraction, with the deformed Lie algebra lying upward.
 The left-pointing arrow symbol over $\mathfrak r_{3,\lambda\neq\pm 1}({\mathbb C})$
 means that it deforms inside the family.

Let us note once again that a non-trivial contraction always induces a
 non-trivial (inverse) jump 
 deformation. The converse is not always true: there are deformations
 which do not admit an inverse contraction. For example, one can never have
 a contraction inside a parameterized family of Lie algebras, but 
 deformations within a family are allowed. Note also that nothing can be contracted to 
 the parameterized family, whereas there are many non-trivial deformations 
 in dimension three to the family $\mathfrak r_{3,\lambda\neq\pm 1}({\mathbb C})$.
 We should emphasize that the irreversibility occurs only when
 we have a family of smooth deformations.
 
The family of Lie algebras $\mathfrak r_{3,\lambda\neq\pm 1}({\mathbb C})$
 has a non-trivial deformation into itself. The two
 Lie algebras $\mathfrak r_{3,1}({\mathbb C})$ and $\mathfrak r_{3,-1}({\mathbb C})$ are special for two reasons. 
 First, $\mathfrak r_{3,1}({\mathbb C})$ can be defor\-med 
 into~$\mathfrak r_{3}({\mathbb C})$, whereas~$\mathfrak r_{3,\lambda\neq 1}({\mathbb C})$ cannot. Moreover
 $\mathfrak r_{3,-1}({\mathbb C})$ can deform into~$\mathfrak{sl}_{2}({\mathbb C})$, 
 whereas $\mathfrak r_{3,\lambda\neq -1}({\mathbb C})$ cannot. Second, $\mathfrak r_{3,1}({\mathbb C})$
 is special because it cannot be contracted to $\mathfrak n_{3}({\mathbb C})$,
 unlike $\mathfrak r_{3,\lambda\neq 1}({\mathbb C})$.

\section[Infinite dimensional Lie algebras]{Inf\/inite dimensional Lie algebras}

 The physical interest of inf\/inite dimensional Lie algebras
 stems mainly from conformal f\/ield
 theory and critical phenomena in two dimensions \cite{cft}.
 Whereas the Witt and Virasoro algeb\-ras describe local invariance
 of conformal f\/ield theories on
 the (zero genus) Riemann sphere, hereafter we shall
 discuss some Lie algebras of Krichever--Novikov
 type which correspond to higher genus. As mentioned
 previously, these Lie algebras are formally
 rigid, yet they can be deformed. These deformations
 are so-called {\it global} and cohomology
 theory then does not lend itself to compute such deformations
 \cite{fialowskimontigny}. These contractions turn out to be
 especially useful since they often lead to new inf\/inite-dimensional
 Lie algebras, as we will demonstrate now.  

\subsection[Witt algebras and Krichever-Novikov algebras]{Witt algebras and Krichever--Novikov algebras}

First, let us consider the Witt algebra ${\mathfrak W}$ with Lie brackets:
\begin{gather}
[l_n,l_m]=(m-n)l_{n+m},\qquad n,m\in {\mathbb Z}.
\label{witt}
\end{gather} 
 The deformation pattern of its only one-dimensional central
 extension, the Virasoro algebra, is quite similar and we will not
 discuss it hereafter.
 Krichever and Novikov introduced new algebras in~\cite{kn}. An interesting aspect of inf\/inite-dimensional Lie
 algebras, which does not occur for the f\/inite-dimensional cases,
 is that Lie algebras of Krichever--Novikov
 type ${\mathfrak{KN}}$ can be interpreted as {\it global}
 deformations of the Witt or Virasoro algebra  \cite{fialowskischli},
 even though the Witt algebra  ${\mathfrak W}$ is formally rigid,
 thus preventing any non-trivial
 {\it formal} deformation.

An example of Krichever--Novikov algebras that can be obtained as a one-parameter
 global deformation of Witt algebra is generated with the following f\/ield basis:
\[
V_{2n}\equiv X(X-\alpha)^n(X+\alpha)^n \frac d{dX},\qquad
 V_{2n+1}\equiv (X-\alpha)^{n+1}(X+\alpha)^{n+1} \frac d{dX},
\]
which satisfy the following Lie brackets:
\begin{gather}
[V_n,V_m]=\left \{\begin{array}{ll}
(m-n)V_{n+m}, & n,m\ {\rm odd},\\[1mm]
(m-n)(V_{n+m}+\alpha^2V_{n+m-2}), & n,m\ {\rm even},\\[1mm]
(m-n)V_{n+m}+(m-n-1)\alpha^2V_{n+m-2},\quad & n\ {\rm odd}, \  m\ {\rm even}.\end{array}
\right.\label{kn1}
\end{gather}

Clearly this can be contracted back to the Witt algebra
 by def\/ining ${\cal U}_\varepsilon$ in
 equations (\ref{wwsplit}) and~(\ref{weimarwoods}) as 
 \begin{gather}
l_n\equiv\varepsilon^n V_n,\qquad {\text {for\ all}} \ n\in{\mathbb Z}. \label{powercontraction}
 \end{gather}
 Then equation (\ref{kn1}) becomes
\begin{gather*}
[l_n,l_m]_\varepsilon=\varepsilon^{n+m}[V_n,V_m]=\left \{\!\begin{array}{ll}
(m-n)l_{n+m},\quad & n,m\ {\rm odd},\\[1mm]
(m-n)(l_{n+m}+\varepsilon^2\alpha^2l_{n+m-2}),\quad & n,m\ {\rm even},\\[1mm]
(m-n)l_{n+m}+(m-n-1)\varepsilon^2\alpha^2l_{n+m-2},\quad & n\ {\rm odd}, \ m\ {\rm even}.\end{array}
\right.\!\!
\end{gather*}
We retrieve the commutation relations, equation (\ref{witt}), of ${\mathfrak W}$ as $\varepsilon$
 approaches zero. Therefore, the operations
 of deformation and contraction are mutually reversible in this case. 

In addition to retrieving the Witt algebra ${\mathfrak{W}}$, equation (\ref{kn1}),  
 one may contract ${\mathfrak{KN}}$ to new Lie algebras.
 For instance, let us def\/ine ${\cal U}_\varepsilon$ as 
 \begin{gather}
 {\cal U}_\varepsilon\equiv \varepsilon^{n_0}{\mathrm {id}}_{\mathfrak{g}_0}+
\varepsilon^{n_1}{\mathrm {id}}_{\mathfrak{g}_1},\label{z2contraction}
 \end{gather}
 where $0$ and $1$ denote the even and odd powers in ${\mathfrak{KN}}$, respectively.
 This choice is quite natural, given the odd versus even splitting in equation (\ref{kn1}).  
 Then the Lie brackets (\ref{kn1}) are modif\/ied to  
 \begin{gather}
 [V_n,V_m]_\varepsilon=\left \{\begin{array}{ll}
 \varepsilon^{2n_1-n_0}(m-n)V_{n+m}, & n,m\ {\rm odd},\\[1mm]
 \varepsilon^{n_0}(m-n)(V_{n+m}+\alpha^2V_{n+m-2}), & n,m\ {\rm even},\\[1mm]
 \varepsilon^{n_0}[(m-n)V_{n+m}+(m-n-1)\alpha^2V_{n+m-2}],\quad
 & n\ {\rm odd}, \ m\ {\rm even}.
\end{array}\right.\label{knmodz2}
 \end{gather}
 Clearly we must have non-negative values of $n_0$ and $2n_1-n_0$.
 This leads to {\it four} separate contracted Lie algebras: 
 (1) we obtain the trivial Abelian Lie algebra
 when these expressions take on positive
 values; (2) another trivial contraction is given by $n_0=n_1=0$;
 then it leaves the commutators of equation (\ref{kn1}) unchanged;
 (3) the In\"on\"u--Wigner contraction,
 given by $n_0=0$ and $n_1=1$, so that
 the contracted commutation relations read:
 \begin{gather*}
 [V_n,V_m]=\left \{\begin{array}{ll}
 0, & n,m\ {\rm odd},\\[1mm]
 (m-n)(V_{n+m}+\alpha^2V_{n+m-2}), & n,m\ {\rm even},\\[1mm]
 (m-n)V_{n+m}+(m-n-1)\alpha^2V_{n+m-2},\quad & n\ {\rm odd}, \ m\ {\rm even}.
 \end{array}\right.%\label{iwkn1}
 \end{gather*}
 and (4) by choosing $n_0>0$ and $2n_1-n_0=0$ in equation (\ref{knmodz2}), we f\/ind 
 \[
 [V_n,V_m]=\left \{\begin{array}{ll}
 (m-n)V_{n+m},\quad & n,m\ {\rm odd},\\[1mm]
 0, & n,m\ {\rm even},\\[1mm]
 0, & n\ {\rm odd}, \ m\ {\rm even}.\end{array}\right.
 \]
The cases (3) and (4) are clearly not isomorphic to the Witt
  algebra, equation (\ref{witt}). 

\subsection[Affine Kac-Moody and Krichever-Novikov algebras]{Af\/f\/ine Kac--Moody and Krichever--Novikov algebras}

Now let us discuss deformations and contractions of ${\mathfrak{KN}}$-type 
 Kac--Moody algebras. Untwisted Kac--Moody algebras 
 $\hat{{\mathfrak{g}}}=( {\mathfrak{g}} \otimes {\mathbb C} [t,t^{-1}] )\oplus {\mathbb C} c$ are def\/ined 
 in terms of a f\/inite simple complex Lie algebra ${\mathfrak{g}}$ together with ${\mathbb C}[t,t^{-1}]$,
 the associative algebra of the Laurent polynomials, and the central extension $c$.
 The Lie brackets may be written as
\begin{gather}
[a\otimes t^m, b\otimes t^n]=[a,b]\otimes t^{m+n}+mcB(a,b)\delta_{m+n,0},
\label{affine}\end{gather}
where $[\;,\;]$ denotes the Lie brackets of the f\/inite Lie algebra ${\mathfrak{g}}$.
 In \cite{majumdar}, Majumdar showed that for In\"on\"u--Wigner contractions,
 af\/f\/inization and contraction are commuting procedures. 

Hereafter we discuss contractions of Krichever--Novikov
 algebras, the latter having been obtained
 by global deformations of af\/f\/ine Kac--Moody algebras.
 In \cite{fialowskischli}, it is shown that the trivially
 extended af\/f\/ine algebras, equation (\ref{affine}) with $k=0$,
 may be deformed to the following $\mathfrak {KN}$ type algebra,
 parameterized over the af\/f\/ine plane ${\mathbb C}^2$, or, described
 algebraically, over the polynomial algebra ${\mathbb C}[e_1,e_2]$:
\begin{gather}
[a\otimes A^n,b\otimes A^m]=\left \{\begin{array}{ll}
[a,b]\otimes A^{n+m},& n\ {\text{or}}\ m\ {\rm even},\\[1mm]
{[a,b]}\otimes A^{n+m}+3e_1\ [a,b]\otimes A^{n+m-2} & \\[1mm]
\quad {} +(e_1-e_2)(2e_1+e_2)[a,b]\otimes A^{n+m-4},\quad & n\
 {\rm and}\ m\ {\rm odd},\end{array}
\right.\label{kn3}
\end{gather} 
where $a$ and $b$ belong to a f\/inite dimensional complex Lie algebra
 ${\mathfrak g}$. The choice $(e_1, e_2)=(0,0)$ simply leads to the original
 af\/f\/ine algebra (although this is {\it not} a contraction process). 

{\samepage As we have done in the previous section, let us f\/irst see that equation 
 (\ref{kn3}) may be contracted back to the original Kac--Moody algebra
 by def\/ining the transformation ${\cal U}_\varepsilon$ in analogy with
 equation (\ref{powercontraction}):
\begin{gather}
a\otimes t^n\equiv\varepsilon^n\ a\otimes A^n,\qquad {\text {for\ all}} \ n\in{\mathbb Z},
 \label{powercontraction2}
 \end{gather}
so that the Lie brackets (\ref{kn3}) become
\[
[a\otimes A^n,b\otimes A^m]_\varepsilon=\left \{\begin{array}{ll}
[a,b]\otimes A^{n+m},\qquad & n\ {\text{or}}\ m\ {\rm even},\\[1mm]
{[a,b]}\otimes A^{n+m}+3e_1\varepsilon^2\ [a,b]\otimes A^{n+m-2} & \\[1mm]
\quad {}+(e_1-e_2)(2e_1+e_2)\varepsilon^4 [a,b]\otimes A^{n+m-4},\quad & n\
 {\rm and}\ m\ {\rm odd}.\end{array}
\right.
\]
This leads to equation (\ref{affine}), with $k=0$, in the limit
 $\varepsilon\rightarrow 0$.} 

Now, let us obtain other Lie algebras
 by using contraction procedures similar to
 what we have done for the Witt algebra. Consider, for instance,
 the splitting of $\mathfrak{KN}$ as in equation (\ref{kn3}) as we have done
 in equation (\ref{z2contraction}). We f\/ind
\begin{gather}
[a\otimes A^n,b\otimes A^m]_\varepsilon=\left \{\begin{array}{ll}
\varepsilon^{n_0}[a,b]\otimes A^{n+m},\qquad & n\ {\text{or}}\ m\ {\rm even},\\[1mm]
\varepsilon^{2n_1-n_0}({[a,b]}\otimes A^{n+m}+3e_1\ [a,b]\otimes A^{n+m-2}\quad & \\[1mm]
\quad {} +(e_1-e_2)(2e_1+e_2)[a,b]\otimes A^{n+m-4}), & n\
 {\rm and}\ m\ {\rm odd}.\end{array}
\right.\label{kn3mod}\hspace{-5mm}
\end{gather}
 The In\"on\"u--Wigner contraction discussed after equation (\ref{knmodz2}),
 for which $n_0=0$ and $n_1=1$, leads to commutation relations
 where the f\/irst line in equation (\ref{kn3mod}), i.e.\
 for $n$ or $m$ even, remains unchanged, whereas the other
 commutators, given in line 2 of the same equation, vanish. 
 If, instead, we take $n_0$ positive and $2n_1-n_0=0$,
 then equation (\ref{kn3mod}) leads to the opposite situation:
 the Lie brackets in the f\/irst line, i.e.\
 for $n$ or $m$ even, will vanish in the limit $\varepsilon\rightarrow 0$,
 whereas the remaining commutators, given in line 2 of 
 equation (\ref{kn3mod}) are left unchanged. 
 Evidently, there are countless possibilities, if we replace equation 
 (\ref{powercontraction2}) with a dif\/ferent splitting of the Lie algebras.

Let us now turn to contractions where the splitting is not done only
 with respect to the degrees of the Laurent polynomials, but  
 within the underlying f\/inite Lie algebra $\mathfrak g$. For the sake of
 illustration, let us consider a ${\mathfrak{KN}}$ algebra based upon the
 Lie algebra ${\mathfrak{sl}_2}({\mathbb C})$, and split it 
 according to the following ${\mathbb Z}_3$-graded structure, inherent
 to the basis of equation (\ref{sl2efh}):  
\[
{\mathfrak{sl}_2}({\mathbb C})={\mathfrak{sl}_2}({\mathbb C})_0+{\mathfrak{sl}_2}({\mathbb C})_1+{\mathfrak{sl}_2}({\mathbb C})_{-1}
 =\{h\}+\{e\}+\{f\}.
\]
We may combine this ${\mathbb Z}_3$-grading with the ${\mathbb Z}_2$-grading provided
 by the even versus odd degrees of the Laurent polynomials to obtain the following 
 ${\mathbb Z}_3\otimes{\mathbb Z}_2$ grading:
\begin{gather*}
{\mathfrak{sl}_2}({\mathbb C})
 =\overbrace{\{h\otimes A^{2n}\}}^{{\mathfrak{sl}_2}({\mathbb C})_{00}}+
\overbrace{\{e\otimes A^{2n}\}}^{{\mathfrak{sl}_2}({\mathbb C})_{10}}+
\overbrace{\{f\otimes A^{2n}\}}^{{\mathfrak{sl}_2}({\mathbb C})_{-10}}+
\overbrace{\{h\otimes A^{2n+1}\}}^{{\mathfrak{sl}_2}({\mathbb C})_{01}}\nonumber\\
\phantom{{\mathfrak{sl}_2}({\mathbb C})=}{}
+ \overbrace{\{e\otimes A^{2n+1}\}}^{{\mathfrak{sl}_2}({\mathbb C})_{11}}+
\overbrace{\{f\otimes A^{2n+1}\}}^{{\mathfrak{sl}_2}({\mathbb C})_{-11}}.
\end{gather*}  
The grading property implies that  
\[
[{\mathfrak g}_\mu,{\mathfrak g}_\nu]_\varepsilon=\varepsilon^{n_\mu+n_\nu-n_{\mu+\nu}}
 {\mathfrak g}_{\mu+\nu},
\]
where $\mu=ab$ is a double index with $a=\{0,1,-1\}\in{\mathbb Z}_3$ and
 $b=\{0,1\}\in{\mathbb Z}_2$. As in equation~(\ref{powerrelation}), the 
 six exponents $n_{00}$,  $n_{10}$, $n_{-10}$, $n_{01}$,
 $n_{11}$, $n_{-11}$ must satisfy $n_\mu+n_\nu-n_{\mu+\nu}\geq 0$. 
 
The commutation relations read explicitly as
\begin{gather}
[h\otimes A^{2n}, e\otimes A^{2m}]_\varepsilon=\varepsilon^{n_{00}} e\otimes A^{2n+2m},\nonumber\\
[h\otimes A^{2n}, e\otimes A^{2m+1}]_\varepsilon=\varepsilon^{n_{00}} e\otimes A^{2n+2m+1},\nonumber\\
[h\otimes A^{2n+1}, e\otimes A^{2m}]_\varepsilon =\varepsilon^{n_{01}+n_{10}-n_{11}} e\otimes A^{2n+2m+1},\nonumber\\
[h\otimes A^{2n+1}, e\otimes A^{2m+1}]_\varepsilon=\varepsilon^{n_{01}+n_{11}-n_{10}} \big(e\otimes A^{2n+2m+2}+
 3e_1 e\otimes A^{2n+2m}\nonumber\\
\phantom{[h\otimes A^{2n+1}, e\otimes A^{2m+1}]_\varepsilon=}{} +(e_1-e_2)(2e_1+e_2) e\otimes A^{2n+2m-2}\big);
\label{hecr}
\\
[h\otimes A^{2n}, f\otimes A^{2m}]_\varepsilon=-\varepsilon^{n_{00}} f\otimes A^{2n+2m},\nonumber\\
[h\otimes A^{2n}, f\otimes A^{2m+1}]_\varepsilon=-\varepsilon^{n_{00}} f\otimes A^{2n+2m+1},\nonumber\\
[h\otimes A^{2n+1}, f\otimes A^{2m}]_\varepsilon=-\varepsilon^{n_{01}+n_{-10}-n_{-11}} f\otimes A^{2n+2m+1},\nonumber\\
[h\otimes A^{2n+1}, f\otimes A^{2m+1}]_\varepsilon=-\varepsilon^{n_{01}+n_{-11}-n_{-10}} \big(f\otimes A^{2n+2m+2}+
 3e_1 f\otimes A^{2n+2m}\nonumber\\
\phantom{[h\otimes A^{2n+1}, f\otimes A^{2m+1}]_\varepsilon=}{} +(e_1-e_2)(2e_1+e_2) f\otimes A^{2n+2m-2}\big);
\label{hfcr}
\\
[e\otimes A^{2n}, f\otimes A^{2m}]_\varepsilon=2 \varepsilon^{n_{10}+n_{-10}-n_{00}} h\otimes A^{2n+2m},\nonumber\\
[e\otimes A^{2n}, f\otimes A^{2m+1}]_\varepsilon=2 \varepsilon^{n_{10}+n_{-11}-n_{01}} h\otimes A^{2n+2m+1},\nonumber\\
[e\otimes A^{2n+1}, f\otimes A^{2m}]_\varepsilon=2 \varepsilon^{n_{11}+n_{-10}-n_{01}} h\otimes A^{2n+2m+1},\nonumber\\
[e\otimes A^{2n+1}, f\otimes A^{2m+1}]_\varepsilon=2 \varepsilon^{n_{11}+n_{-11}-n_{00}} \big(h\otimes A^{2n+2m+2}+
 3e_1 h\otimes A^{2n+2m}\nonumber\\
\phantom{[e\otimes A^{2n+1}, f\otimes A^{2m+1}]_\varepsilon=}{} +(e_1-e_2)(2e_1+e_2) h\otimes A^{2n+2m-2}\big).
\label{efcr}
\end{gather}

We just illustrate a few original algebraic objects that can be obtained
 by contractions. It is not our purpose to f\/ind all the solutions for the
 $n$'s. The f\/irst two lines of equations (\ref{hecr}) and~(\ref{hfcr}) imply
 that 
\[ n_{00}\geq 0.\]
If we take $n_{00}=0$ and $n_{01}=0$, then lines 3 and 4 of equation (\ref{hecr}) imply
 that $n_{11}=n_{10}$, whereas lines 3 and 4 of equation (\ref{hfcr}) lead to
 $n_{-10}=n_{-11}$. From the factors in equation (\ref{efcr}), we f\/ind
 $n_{10}+n_{-10}\geq 0$. All these result in the fact that the commutators
 in equations (\ref{hecr}) and (\ref{hfcr}) all remain unchanged, whereas the
 Lie brackets in equation (\ref{efcr}) either all remain
 unchanged after the contraction, or they all vanish. Note that the 
 latter contracted algebra can also be obtained by f\/irst contracting
 ${\mathfrak{sl}_2}({\mathbb C})$ to ${\mathfrak{r}_{3,-1}}({\mathbb C})$ and then af\/f\/inizing
 it {\it \`a la} Krichever--Novikov. In other words, it follows from the
 Krichever--Novikov construction that the ${\mathfrak{KN}}$ algebra
 obtained directly from ${\mathfrak{r}_{3,-1}}({\mathbb C})$ is the same as
 is we construct the ${\mathfrak{KN}}$ of ${\mathfrak{sl}_2}({\mathbb C})$, and then
 form its contraction, in a way analogous to the contraction from
 ${\mathfrak{sl}_2}({\mathbb C})$ to ${\mathfrak{r}_{3,-1}}({\mathbb C})$. 

Clearly, more complicated contracted algebras, where the commutators
 involving dif\/ferent powers contract dif\/ferently, can be obtained.
 Consider, once again, the case $n_{00}=0$, but $n_{01}>0$. Then,
 the exponents $n_{01}+n_{10}-n_{11}$, $n_{01}+n_{11}-n_{10}$,
 $n_{01}+n_{-10}-n_{-11}$ or $n_{01}+n_{-11}-n_{-10}$ cannot be
 all equal to zero simultaneously. Consider the case where these
 exponents are all strictly positive. 
 This means that lines 1 and 2 of equations (\ref{hecr}) and
 (\ref{hfcr}) remain unchanged under these contractions, whereas
 lines 3 and 4 of the same equations will vanish in the
 limit $\varepsilon\rightarrow 0$. We may choose to preserve lines
 1 and 4 of equation (\ref{efcr}) by taking $n_{10}+n_{-10}=0$ 
 and $n_{11}+n_{-11}=0$, respectively. For the sake of
 illustration, let us choose
\[
n_{00}=0,\qquad n_{01}=n_{10}=n_{11}=1,\qquad
 n_{-10}=n_{-11},
\]
which satisfy all these conditions. Then equations (\ref{hecr}),
 (\ref{hfcr}) and (\ref{efcr}) read
\begin{gather*}
[h\otimes A^{2n}, e\otimes A^{2m}]_\varepsilon=\varepsilon^{0} e\otimes A^{2n+2m},\\
[h\otimes A^{2n}, e\otimes A^{2m+1}]_\varepsilon=\varepsilon^{0} e\otimes A^{2n+2m+1},\\
[h\otimes A^{2n+1}, e\otimes A^{2m}]_\varepsilon=\varepsilon  e\otimes A^{2n+2m+1},\\
[h\otimes A^{2n+1}, e\otimes A^{2m+1}]_\varepsilon=\varepsilon  \big(e\otimes A^{2n+2m+2}+
 3e_1  e\otimes A^{2n+2m}\\
 \phantom{[h\otimes A^{2n+1}, e\otimes A^{2m+1}]_\varepsilon=}{} +(e_1-e_2)(2e_1+e_2)  e\otimes A^{2n+2m-2}\big);
\\
[h\otimes A^{2n}, f\otimes A^{2m}]_\varepsilon=-\varepsilon^{0}  f\otimes A^{2n+2m},\\
[h\otimes A^{2n}, f\otimes A^{2m+1}]_\varepsilon=-\varepsilon^{0}  f\otimes A^{2n+2m+1},\\
[h\otimes A^{2n+1}, f\otimes A^{2m}]_\varepsilon=-\varepsilon  f\otimes A^{2n+2m+1},\\
[h\otimes A^{2n+1}, f\otimes A^{2m+1}]_\varepsilon=-\varepsilon  \big(f\otimes A^{2n+2m+2}+
 3e_1  f\otimes A^{2n+2m}\\
\phantom{[h\otimes A^{2n+1}, f\otimes A^{2m+1}]_\varepsilon=}{} +(e_1-e_2)(2e_1+e_2)  f\otimes A^{2n+2m-2}\big);
\\
[e\otimes A^{2n}, f\otimes A^{2m}]_\varepsilon=2  \varepsilon^0  h\otimes A^{2n+2m},\\
[e\otimes A^{2n}, f\otimes A^{2m+1}]_\varepsilon=2  \varepsilon  h\otimes A^{2n+2m+1},\\
[e\otimes A^{2n+1}, f\otimes A^{2m}]_\varepsilon=2  \varepsilon  h\otimes A^{2n+2m+1},\\
[e\otimes A^{2n+1}, f\otimes A^{2m+1}]_\varepsilon=2  \varepsilon^{0}  \big(h\otimes A^{2n+2m+2}+
 3e_1  h\otimes A^{2n+2m}\\
\phantom{[e\otimes A^{2n+1}, f\otimes A^{2m+1}]_\varepsilon=}{} +(e_1-e_2)(2e_1+e_2)  h\otimes A^{2n+2m-2}\big).
\end{gather*}
In the limit $\varepsilon\rightarrow 0$, they become
\begin{gather*}
[h\otimes A^{2n}, e\otimes A^{2m}]'= e\otimes A^{2n+2m},\\
[h\otimes A^{2n}, e\otimes A^{2m+1}]'= e\otimes A^{2n+2m+1},\\
[h\otimes A^{2n+1}, e\otimes A^{2m}]'=0,\\
[h\otimes A^{2n+1}, e\otimes A^{2m+1}]'=0;\\
[h\otimes A^{2n}, f\otimes A^{2m}]'=-f\otimes A^{2n+2m},\\
[h\otimes A^{2n}, f\otimes A^{2m+1}]'=-f\otimes A^{2n+2m+1},\\
[h\otimes A^{2n+1}, f\otimes A^{2m}]'=0,\\
[h\otimes A^{2n+1}, f\otimes A^{2m+1}]'=0;\\
[e\otimes A^{2n}, f\otimes A^{2m}]'=2h\otimes A^{2n+2m},\\
[e\otimes A^{2n}, f\otimes A^{2m+1}]'=0,\\
[e\otimes A^{2n+1}, f\otimes A^{2m}]'=0,\\
[e\otimes A^{2n+1}, f\otimes A^{2m+1}]'=2  \big(h\otimes A^{2n+2m+2}+
 3e_1  h\otimes A^{2n+2m}\\
\phantom{[e\otimes A^{2n+1}, f\otimes A^{2m+1}]'=}{}
 +(e_1-e_2)(2e_1+e_2)  h\otimes A^{2n+2m-2}\big).
\end{gather*}
This contraction is interesting because it cannot be seen as
 the result of an af\/f\/inization of a contraction of 
 ${\mathfrak{sl}}_2({\mathbb C})$. Indeed, it is not an af\/f\/inization of
 ${\mathfrak{r}}_{3,-1}({\mathbb C})$, in which case all the commutation
 relations of type $[e,f]$, in equation (\ref{efcr}), would vanish, unlike
 the commutation relations above. This is similar for 
 ${\mathfrak{n}}_3({\mathbb C})$, for which the Lie brackets of types
 $[h,e]$ and $[h,f]$ (equations (\ref{hecr}), (\ref{hfcr})) do vanish,
 unlike the Lie brackets obtained above.
 Also, we cannot see the contractions above as involving only the
 even versus odd powers of the Laurent polynomials because the
 grading involves also the underlying f\/inite-dimensional Lie
 algebra. Indeed we f\/ind in the contracted Lie algebra that the
 Lie brackets are not uniquely determined by the parity of their
 elements; for instance, some commutators of elements both even
 do commute, whereas other such commutators do not commute.   

Clearly, numerous other limits can be obtained, but our purpose here
 was just to demonstrate how a combination of deformations
 and contractions can lead to new Lie algebras. 
 
\subsection*{Acknowledgements}

Partial support was provided by a grant from NATO Hungary.
 AF is grateful to OTKA (Hungary) for grants T043641 and T043034.
 MdM acknowledges partial f\/inancial support from the Natural Sciences and
 Engineering Research Council of Canada, and the organizers of the
 Conference Symmetry in Nonlinear Mathematical Physics, held in
 Kyiv in June 2005. We thank both referees for their remarks and
 useful suggestions.

\LastPageEnding


\begin{thebibliography}{99}
\footnotesize

\bibitem{fialowskimontigny} Fialowski A., de Montigny M., Contractions
 and deformations of Lie algebras, {\it J. Phys. A: Math. Gen.}, 2005, V.38, 6335--6349.

\bibitem{gerstenhaber} Gerstenhaber M., On the deformation of rings and algebras,
{\it Ann. Math.}, 1964, V.79, 59--103;\\
Gerstenhaber M., On the deformation of rings and algebras II, {\it Ann. Math.}, 1966, V.84, 1--19;\\ 
Gerstenhaber M., On the deformation of rings and algebras III, {\it Ann. Math.}, 1968, V.88, 1--34;\\ 
Gerstenhaber M., On the deformation of rings and algebras IV, 1974, V.99, 257--276;\\
Nijenhuis A., Richardson R.W., Deformations of Lie algebra structures, 
{\it J. Math. Mech.}, 1967, V.17, 89--105;\\
Fialowski A.,  Deformations of Lie algebras, {\it  Math. USSR Sbornik}, 1986, V.55, 467--472;\\
Fialowski A., An example of formal deformations of
 Lie algebras, in Proceedings of NATO Conference on
 Deformation Theory of Algebras and Applications, Editors M.~Hazawinkel and
 M.~Gerstenhaber, Dordrecht, Kluwer, 1988, 375--401;\\
 Fialowski A., Fuchs D., Construction of miniversal deformations
 of Lie algebras, {\it J. Funct. Anal.}, 1999, V.161, 76--110, math.RT/0006117.

\bibitem{fialowskischli} Fialowski A., Schlichenmaier M., Global deformations
 of the Witt algebra of Krichever--Novikov type, {\it Commun. Contemp. Math.}, 2003, V.5, 921--945, math.QA/0206114;\\
Fialowski A., Schlichenmaier M., Global geometric deformations of current
 algebras as Krichever--Novikov type algebras, {\it Comm. Math. Phys.}, 2005, V.260, 579--612,
 math.QA/0412113. 

\bibitem{iw} In\"on\"u E., Wigner E.P., On the contraction of groups and their
 representations, {\it Proc. Nat. Acad. Sci. U.S.A.}, 1953, V.39, 510--524;\\
Saletan E., Contraction of Lie groups, {\it J. Math. Phys.}, 1961, V.2, 1--21;\\
Gilmore R., Lie groups, Lie algebras, and some of their
applications, New York, Wiley, 1974, Chapter 10;\\
Talman J.D., Special functions: a group theoretic approach,
New York, Benjamin,  1968.

\bibitem{vinberg} Onishchik A.L., Vinberg E.B.,  Lie groups and
 Lie algebras, {\it Enclycopaedia of Mathematical Sciences}, Vol.~41,
Berlin, Springer, 1991, Chapter~7.

\bibitem{weimar2000} Weimar-Woods E., Contractions, generalized
 In\"on\"u--Wigner contractions and deformations of f\/inite-di\-men\-sional
 Lie algebras, {\it Rev. Math. Phys.}, 2000, V.12, 1505--1529;\\
L\'evy-Nahas M., Deformation and contraction
 of Lie algebras, {\it J. Math. Phys.}, 1967, V.8, 1211--1222;\\
L\~ohmus J., Tammelo R., Contractions and deformations of space-time
 algebras I. General theory and kinematical algebras, {\it Hadronic J.}, 1997, V.20, 361--416;\\
Fialowski A., O'Halloran J., A comparison of deformations
 and orbit closure, {\it Comm. Algebra}, 1990, V.18, 4121--4140.

\bibitem{conatser} Conatser C.W., Contractions of the low-dimensional real
 Lie algebras, {\it J. Math. Phys.}, 1972, V.13, 196--203.

\bibitem{burdesteinhoff} Burde D., Steinhof\/f C., Classif\/ication of orbit closures of
 4-dimensional complex Lie algebras, {\it J. Algebra}, 1999, V.214, 729--739.
 
 
\bibitem{fialowski1990} Fialowski A., Deformations of some
 inf\/inite-dimensional Lie algebras, {\it J. Math. Phys.}, 1990, V.31, 1340--1343.

\bibitem{formalrigid}  Lecomte P.B.A., Roger C.,
Rigidity of current Lie algebras of complex simple type, {\it J. London Math. Soc.~(2)},
1988, V.37, 232--240.

\bibitem{fialowskipenkava} Fialowski A., Penkava M.,
 Versal deformations of three-dimensional Lie algebras as
 $L_\infty$ algebras, {\it Commun. Contemp. Math.}, 2005, V.7, 145--165, math.RT/0303346.


\bibitem{cft} Goddard P., Olive D., Kac--Moody and Virasoro algebras,
New York, World Scientif\/ic, 1988;\\
Di Francesco P., Mathieu P., S\'en\'echal D., Conformal f\/ield theory,
New York, Springer, 1997; \\
Tsvelik A.M., Quantum f\/ield theory in condensed matter physics,
Cambridge Univ. Press, 2003.

\bibitem{kn} Krichever I.M., Novikov S.P., Algebras of Virasoro type,
  Riemann surfaces and structures of the theory of solitons, {\it Funct. Anal. Appl.}, 1987, V.21, 126--142;\\ 
Krichever I.M., Novikov S.P., Virasoro-type algebras, Riemann surfaces
  and strings in Minkowski space, {\it Funct. Anal. Appl.}, 1987, V.21, 294--307;\\
Krichever I.M., Novikov S.P., Algebras of Virasoro type,
  energy-momentum tensor and decomposition ope\-rators on Riemann surfaces, {\it Funct. Anal. Appl.},
1989, V.23, 19--33.

\bibitem{majumdar} Majumdar P., In\"on\"u--Wigner contraction of Kac--Moody algebras,
{\it J. Math. Phys.}, 1993, V.34, 2059--2065, hep-th/9207057;\\
Olive D.I., Rabinovici E., Schwimmer A., A class of string backgrounds as
 a semiclassical limit of WZW models, {\it Phys. Lett. B}, 1994, V.321, 361--364, hep-th/9311081.
\end{thebibliography}
\end{document}